\documentclass[12pt]{amsart}%{article}%
\usepackage{amsmath,amstext,amsfonts,amsthm,amssymb}
\usepackage{eucal}
\usepackage{graphicx}
\renewcommand{\subsubsection}[1]{\addtocounter{subsubsection}{1}
{\ \\[3pt]\bf \thesubsubsection. \  #1} }
\swapnumbers

{  \theoremstyle{definition}

}
\newcommand{\End}{\operatorname{End}}
\newcommand{\Id}{\operatorname{Id}}

\newcommand{\Hom}{\operatorname{Hom}}

\newcommand{\lra}{\longrightarrow}
\newcommand{\ra}{\rightarrow}

\newcommand{\iso}{\overset\sim\longrightarrow}

\newcommand{\dpar}{\partial}

\newcommand{\bea}{\begin{eqnarray*}}
\newcommand{\eea}{\end{eqnarray*}}
\newcommand{\bean}{\begin{eqnarray}}
\newcommand{\eean}{\end{eqnarray}}

\newcommand{\fg}{\mathfrak g}
\newcommand{\fgl}{\mathfrak{gl}}
\newcommand{\fG}{\mathfrak G}

\newcommand{\fk}{\mathfrak k}
\newcommand{\fK}{\mathfrak K}

\newcommand{\fsl}{\mathfrak{sl}}

\newcommand{\CB}{\mathcal{B}}
\newcommand{\CC}{\mathcal{C}}

\newcommand{\CK}{\mathcal{K}}
\newcommand{\CL}{\mathcal{L}}

\newcommand{\CY}{\mathcal{Y}}

\newcommand{\BC}{\mathbb{C}}

\newcommand{\BQ}{\mathbb{Q}}
\newcommand{\BR}{\mathbb{R}}
\newcommand{\BZ}{\mathbb{Z}}

\begin{document}

%  WRITE HERE

\centerline{\bf INVARIANT TRIPLE FUNCTIONALS}

\bigskip\bigskip 

\centerline{\bf OVER $U_q\frak{sl}_2$} 

\bigskip\bigskip

\centerline{Bui Van Binh and Vadim Schechtman}

\bigskip\bigskip

%\vspace{1cm}

%\input loke-q-0

\centerline{\bf Introduction}

\bigskip\bigskip

Before describing the contents of this note let us discuss  
some motivation and questions behind it. 

The fact that an irreducible finite dimensional representation 
$V(\lambda_1)$ of highest weight $\lambda_1$ of the Lie algebra $\fg = \fsl_2(\BC)$ 
occurs with multiplicity at most $1$ in a  tensor product 
$V(\lambda_2)\otimes V(\lambda_3)$ is easy and classical. Since these representations 
are isomorphic to their duals, the same thing may be expressed 
by saying that the dimension of the space of $\fg$-invariant 
functionals 
$$
\dim \Hom_\fg(V(\lambda_1)\otimes V(\lambda_2)\otimes V(\lambda_3),\BC) \leq 1
\eqno{(0.1)}
$$ 
The multiplicity one statements like (0.1) hold true as well if $V(\lambda_i)$ are irreducible 
infinite dimensional representations of real, complex and $p$-adic Lie groups or Lie algebras 
close to $GL_2$  (their proof being usually more difficult). 

As an example, such a statement for the group $G = PGL_2(\BR)$ and the representations 
of the principal series is applied in [BR]. 
In that case a representation $V(\lambda)$ may be realized (before the Hilbert completion) 
in the space of smooth functions on the unit circle $f:\ S^1 \lra \BC$, and the tensor product 
$V(\lambda_1)\otimes V(\lambda_2)\otimes V(\lambda_3)$ --- in the space of functions of three variables 
$f:\ (S^1)^3\lra \BC$. An explicit linear functional
$$
\ell_{\lambda_1,\lambda_2,\lambda_3}:\ 
V(\lambda_1)\otimes V(\lambda_2)\otimes V(\lambda_3) \lra \BC
$$
may be defined in the form of an integral
$$
\ell_{\lambda_1,\lambda_2,\lambda_3}(f) = \int_{(S^1)^3} f(\theta_1,\theta_2,\theta_3)
\fK_{\lambda_1,\lambda_2,\lambda_3}(\theta_1,\theta_2,\theta_3)
d\theta_1d\theta_2d\theta_3
\eqno{(0.2)}
$$ 
against some naturally defined $G$-invariant kernel 
$\fK_{\lambda_1,\lambda_2,\lambda_3}$, 
cf. [BR], 5.1.1, [Ok] (0.10), (0.12). On the other hand 
our triple product contains a distinguished spherical 
(i.e. $PO(2)^3$-invariant) vector $v_{\lambda_1,\lambda_2,\lambda_3}$, 
the constant function $1$. 

The value 
$$
\ell_{\lambda_1,\lambda_2,\lambda_3}(v_{\lambda_1,\lambda_2,\lambda_3}) = 
\int_{(S^1)^3} 
\fK_{\lambda_1,\lambda_2,\lambda_3}(\theta_1,\theta_2,\theta_3)
d\theta_1d\theta_2d\theta_3
\eqno{(0.3)}
$$  
is 
equal to certain quotient of products of Gamma values. Its asymptotics with respect to $\lambda_i$ 
(which follows from the Stirling formula) is one of the ingredients used in [BR] for an estimation of Fourier coefficients 
of automorphic triple products. 

In the paper [BS] we have calculated the integrals similar to (0.3)  
corresponding to complex and $p$-adic groups $PGL_2(\BC)$, $PGL_2(\BQ_p)$, and also an analogous 
{\it $q$-deformed} integral which has the form
$$
\int_{(S^1)^3} 
\fK_{\lambda_1,\lambda_2,\lambda_3;q}(\theta_1,\theta_2,\theta_3)
d\theta_1d\theta_2d\theta_3
\eqno{(0.4)}
$$
where $\fK_{\lambda_1,\lambda_2,\lambda_3;q}$ is a certain 
$q$-deformation of the kernel 
$\fK_{\lambda_1,\lambda_2,\lambda_3}$. These integrals  
are expressed in terms of the complex, $p$-adic and $q$-deformed  versions 
of $\Gamma$-functions respectively. One could expect that it is possible 
to find representations $V_q(\lambda)$ of the $q$-deformed algebra 
$U_q\fgl_2$ in the space of functions on $S^1$, so that the 
$q$-deformed kernel $\fK_{\lambda_1,\lambda_2,\lambda_3;q}$ 
will be a $U_q\fgl_2$-invariant element of the triple 
product 
$V_q(\lambda_1)\otimes V_q(\lambda_2)\otimes V_q(\lambda_3)$. 

We do not pursue this direction here, but we prove some multiplicity one statement like 
(0.1) over the quantum group. 
Our starting point was a theorem of Hung Yean Loke [L] 
who proves in particular that (0.1) holds true if $\fg = \fgl_2(\BR)$ and $V(\lambda_i)$ are  
irreducible representations of the (infinitesimal) principal series 
defined by Jacquet-Langlands, cf. [JL], Ch. I, \S 5. 
The space of such a representation is much smaller than 
the spaces of smooth functions above, it is rather a "discrete analog" of it, and the structures that 
appear are quite similar. 

A base $\{e_q,\ q\in Q\}$ of $V(\lambda)$ is enumerated by a 
set $Q$ which may be identified with coroot lattice of $\fg$ 
(or with $\BZ$). Thus elements of $V(\lambda)$ are 
finite linear combinations 
$$
\sum a(q)e_q
$$
which we can consider as functions $a:\ Q\lra \BC$ which are 
compactly supported, i.e. all but finitely many $a(q)$ are zeros. 
The Lie algebra $\fg$ acts on these functions by 
difference operators (depending on $\lambda\in \BC$\footnote{there is also a discrete parameter $\epsilon$ around which is not important 
in our discussion, so we forget about it in this Introduction}) 
of order $\leq 1$ 
(to avoid the confusion, $V(\lambda)$ is {\it not} 
a highest weight module). Thus,  elements 
of a triple product $V(\lambda_1)\otimes V(\lambda_2)\otimes V(\lambda_3)$ are compactly supported functions 
$a:\ Q^3\lra \BC$. Similarly, a trilinear functional        
$$
\ell:\ V(\lambda_1)\otimes V(\lambda_2)\otimes V(\lambda_3) \lra \BC
$$
is uniquely determined by its action on the basis elements. 
If we denote 
$$
k(q_1,q_2,q_3) = \ell(e_{q_1}\otimes e_{q_2}\otimes e_{q_3}),
$$
we get a function $\CK:\ \BQ^3\lra\BC$ (an arbitrary, not necessarily compactly supported one). 
The value of $\ell$ is given by
$$
\ell(a) = \sum_{(q_1,q_2,q_3)\in Q^3} a(q_1,q_2,q_3)k(q_1,q_2,q_3);
$$
this formula is similar to (0.2). 

The functional is $\fg$-invariant iff the corresponding 
$\CK$ satisfies a simple system of difference equations. 
The result of [L] says that the space of such functions $\CK$  
is one-dimensional. It would be interesting to find 
a nice explicit formula for a solution.  
   
In \S 2 of the present note we define 
principal series representations over the quantum group $U_q\fsl_2$ 
which are $q$-deformations of the Jacquet - Langlands modules. 
Then we  
define natural intertwining ("reflection") operators between them (cf. 2.2.2) and finally prove 
for them an analog of (0.1), cf. Thm. 2.3 for the precise formulation; this is our main result.  
The proof is a $q$-deformation of the argument from [L]. 

In \S 1 we recall the definitions from [JL] and the original argument of [L] and present some comments on it, 
cf. 1.4, in the spirit of I.M.Gelfand's philosophy considering 
the Clebsch-Gordan coefficients as discrete orthogonal polynomials, cf. [NSU].  

We thank F.Malikov who has drawn our attention to a very interesting paper [FM].

\bigskip\bigskip

\centerline{\bf \S 1. Invariant triple functionals over $\frak{sl}_2$} 

\bigskip\bigskip

\bigskip

{\bf 1.1. Principal series.} First we recall the classical definition of the principal series following Jacquet - Langlands. Another definition  
of these modules may be found in [FM]. 

Let $\fg = \frak{sl}_2$ and $E, F, H$ be the standard base of $\fg$. 

Let $s\in \BC$, $\epsilon\in \{0, 1\}$. 
Following [JL], \S 5 and [L] 2.2  consider the following representation $M(s,\epsilon)$ of $\fg$ (for a motivation of the definition 
see 1.6 below.) 

The underlying vector space of $M(s,\epsilon)$ has a $\BC$-base $\{v_n\}_{n\in \epsilon + 2\BZ}$. We denote $M_n = \BC\cdot v_n$, 
so $M(s,\epsilon) = \oplus_{n\in\BZ} M_n$ where we set $M_n = 0$ if $n\notin \epsilon + 2\BZ$.
  
The action of $\fg$ is given by
$$
Hv_n = n v_n,\ 
\eqno{(1.1.1)}
$$
$$
Ev_n = \frac{1}{2}(s + n + 1)v_{n+2},\ Fv_n = \frac{1}{2}(s - n + 1)v_{n-2}
\eqno{(1.1.2)}
$$
Thus, 
$$
[H,E] = 2E,\ [H,F] = - 2F,\ [E, F] = H
$$

{\bf 1.1.1. Lemma.} {\it If $s - \epsilon\notin 2\BZ + 1$ then 
$M(s,\epsilon)$ is an irreducible $\fg$-module.}

{\bf Proof.} Let $W\subset M:= M(s,\epsilon)$ be a $\fg$-submodule, $W\neq 0$. Since $W$ is $H$-invariant, $W = \oplus_n W\cap M_n$. 
Thus there exists $x\in W\cap M_n, x\neq 0$. Due to the hypothesis $F^mx\neq 0$ and $E^mx\neq 0$ for all $m\in\BZ$, whence 
$W = M$. $\square$

{\bf 1.1.2. The reflection operator.} Cf. [JL], between 5.11 and 5.12. Consider two modules $M(\pm s, \epsilon)$. A linear map 
$$
f:\ M(s,\epsilon) \lra M(-s,\epsilon)
$$
is $\fg$-equivariant iff it respects the gradings (since it commutes with $H$), say $f(v_n) = f_n v'_n$, and the numbers $f_n$ 
satisfy two relations
$$
(s + n + 1)f_{n+2} = (-s + n + 1)f_{n},
\eqno{(1.1.3a)}
$$
(commutation with $E$) and 
$$
(s - n + 1)f_{n-2} = (-s - n + 1)f_{n}
\eqno{(1.1.3b)}
$$
(commutation with $F$). In fact these equations are equivalent: for example 
$(1.1.3b)$ is the same as $(1.1.3a)$ with $n$ replaced by $n-2$, multiplied by $-1$. 

These relations are satisfied if 
$$
f_n = \frac{\Gamma((-s+n+1)/2)}{\Gamma((s+n+1)/2)},
$$
We shall denote the corresponding intertwining operator by 
$$
R(s):\ M(s,\epsilon) \iso M(-s,\epsilon)
$$
In fact, these are the only possible intertwiners between 
different modules of principal series. 

One has
$$
R(-s)R(s) = \Id_{M(s)}.
$$

{\bf 1.2. Theorem}, cf. [L], Thm 1.2 (1). {\it Consider three $\fg$-modules 
$M^i = M(s_i, \epsilon_i)$, $i = 1, 2, 3$. Suppose that $s_i - \epsilon_i\notin 1 + 2\BZ$. 
There exists a unique, up to a multiplicative constant, function 
$$
f:\ M:= M^1\otimes M^2\otimes M^3 \lra \BC
$$
such that
$$
f(X m) = 0,\ X\in\fg,\ m\in M
\eqno{(1.2.1)}
$$
and
$$
f(\omega m) = f(m)
\eqno{(1.2.2)}
$$
where $\omega:\ M\iso M$ is an automorphism defined by 
$$
\omega(v_n\otimes v_m\otimes v_k) = v_{-n}\otimes v_{-m}\otimes v_{-k}.
$$}

{\bf 1.3. Proof} (sketch). The condition (1.2.1) for $X = H$ implies that 
$f(v_n\otimes v_m\otimes v_k) = 0$ unless $n + m + k = 0$. 

Let us denote
$$
a(n, m) = f(v_n\otimes v_m\otimes v_{-n-m})
$$

The conditions (1.2.1) and (1.2.2) are equivalent to a system of $3$ equations on the  function $a(n,m)$: 
$$
a(n, m) = a(-n, -m),
\eqno{(1.3.0)}
$$
$$
(s_1 + n + 1)a(n+2,m) + (s_2 + m + 1)a(n,m+2) 
+ (s_3 - n - m - 1)a(n,m) = 0
\eqno{(1.3.1)}
$$
and 
$$
(s_1 - n + 1)a(n-2,m) + (s_2 - m + 1)a(n,m-2) + (s_3 + n + m - 1)a(n,m) = 0
\eqno{(1.3.2)}
$$
One has to show that these equations admit a unique, up to scalar, solution.

By considering the "bonbon" configuration 
$$
B = \{(n,m),(n-2,m),(n,m+2),(n-2,m+2),(n+2,m),(n,m-2),(n+2),(n+2)\}
$$
one shows\footnote{there is a misprint in [L]: in formula (5) on p. 124 one should interchange 
$c$ with $d$, and in formula (6) --- $d$ with $e$.} that (1.3.1-2) imply an equation
$$
(s_1 - n + 1)(s_2 + m + 1)a(n-2,m+2) - 
(s_3^2 - s_1^2 - s_2^2 - 2nm + 1) a(n,m) +   
$$
$$
(s_1 + n + 1)(s_2 - m + 1)a(n+2,m-2) = 0
\eqno{(1.3.3)} 
$$
After that it is almost evident that a solution of (1.3.1-2) is uniquely defined by its two values 
on a diagonal, like $a(n,m), a(n-2,m+2)$. The parity condition (1.3.0) implies 
that the space of solutions has dimension $\leq 1$.

The non-trivial part is a proof of the {\it existence} of a solution. It is 
a direct computation.  
Cf. the argument for the 
$q$-deformed case in \S 2 below. $\square$

{\bf 1.4. Difference equations on the root lattice of type $A_2$.} Let $X$ denote the lattice $\{(n_1,n_2)\in\BZ^2|\ n_i - \epsilon_i\in 2\BZ\}$. (Note that initially 
it comes in the above proof as a lattice 
$$
\{(n_1,n_2,n_3)\in\BZ^3|\ n_i - \epsilon_i\in 2\BZ,\ \sum n_i = 0\}
$$
and resembles the root lattice of the root system of type $A_2$.) 

Consider the space of maps of sets  
$Y = \{a:\ X\lra \BC\}$; $Y$ is a complex vector space. Define two linear operators 
$L_\pm\in \End Y$ by 
$$
L_+a(n,m) = 
$$
$$
(p + n)a(n+2,m) + (s+n)a(n,m+2) + (r-n-m)a(n,m),
\eqno{(1.4.1a)}
$$ 
$$
L_-a(n,m) = 
$$
$$
(p - n)a(n-2,m) + (s-n)a(n,m-2) + (r+n+m)a(n,m),
\eqno{(1.4.1b)}
$$
where
$p = s_1 + 1, s = s_2 + 1, r = s_3 - 1$. 

One can rewrite the equations (1.3.1-2) in the form
$$
L_+ a = 0,\ L_-a = 0
\eqno{(1.4.2)}
$$
{\bf 1.4.1. Lemma.} $[L_+, L_-] = 2(L_+ - L_-)$. $\square$

It follows that $L_+$ and $L_-$ span a $2$-dimensional Lie algebra isomorphic to a Borel subalgebra 
of $\fsl_2$. 

Following [NSU], Ch. II, \S 1, introduce   
the forward and backward difference ("discrete derivatives") operators acting on functions $f(n)$ of an integer argument:
$$
\Delta  f(n) = f(n+2) - f(n),\ \nabla f(n) = f(n) - f(n-2)
$$
These operators give rise to "discrete partial derivatives" acting on 
the space of functions of two variables $a(n,m)$ as above. We denote by subscripts $_n$ or $_m$ the operators acting on the first 
(resp. second) argument, for example 
$$
\nabla_na(n,m) = a(n,m) - a(n-2,m), 
$$
etc. 
Then the equations (1.4.2) rewrite as follows: 
$$
((n+p)\Delta_n + (m+s)\Delta_m)a = -(p+s+r)a
\eqno{(1.4.3a)}
$$
$$
((n-p)\nabla_n + (m-s)\nabla_m)a = -(p+s+r)a
\eqno{(1.4.3b)}
$$
These equations are similar to [NSU], Ch. IV, \S 2, (30). 

Let us fix $k$ and consider the functions $b(n) := a(n,k-n)$. The equation (1.3.3) 
is a second order equation satisfied by these functions may be written as
$$
\bigl\{(n+p)(n+s-k)\nabla\Delta + 2(pn + sn - pk)\nabla - r(r-2)\bigr\}b = 0
\eqno{(1.4.4)}
$$ 
It is a "difference equation of hypergeometric type" in the terminology of 
[NSU], Ch. II, \S 1. Their solutions can be called "Hahn functions".

{\bf 1.5. Analogous differential equations.}  
It is instructive to consider the continuous analogs of the previous operators. 

Let us consider the following operators on the space $\CY$ of differentiable functions 
$a(x,y):\ \BR^2 \lra \BC$ which is a continuous analog of the space $Y$:
$$
\CL_+ = (p + x)\dpar_x + (s + y)\dpar_y + p + s +  r
$$
$$
\CL_- = (-p + x)\dpar_x + (-s + y)\dpar_y + p + s +  r
$$

{\bf 1.5.1. Lemma.} $[\CL_+, \CL_-] = \CL_+ - \CL_-$. $\square$

The analog of (1.4.4) is a hypergeometric equation 
$$
(x+p)(x+q-k)b''(x) + 2((p+s)x - pk)b'(x) - r(r-2)b(x) = 0
\eqno{(1.5.1)}
$$
where $b(x) = a(x,k-x)$.

\bigskip

{\bf 1.6. Motivation: Jacquet - Langlands principal series over $GL_2(\BR)$.} Cf. [GL], Ch. I, \S 5. 
Recall that a {\it quasicharacter} of the group 
$\BR^*$ is a continuous homomorphism $\mu:\ \BR^* \lra \BC^*$. All such homomorphisms have the form 
$$
\mu(x) = \mu_{s,m}(x) = |x|^s(x/|x|)^m
$$
where $s\in \BC,\ m \in\{0, 1\}$. Let $\mu_i = \mu_{s_i,m_i}$, $i = 1, 2$, be two such quasicharacters. Let 
$\CB'(\mu_1,\mu_2)$ denote the space of $C^\infty$-functions $f: G := GL_2(\BR)\lra \BC$ such that
$$
f\biggl(\left(\begin{matrix} a & c\\ 0 & b\end{matrix}\right)g\biggr) = \mu_1(a)\mu_2(b)(|a/b|)^{1/2}f(g)
$$
for all $g\in G, a, b\in \BR^*, c\in \BR$. $G$ acts on $\CB'(\mu_1,\mu_2)$ in the obvious way. 

Set $s = s_1 - s_2$ and $m = |m_1 - m_2|$. For any $n\in\BZ$ such that $n - m\in 2\BZ$ define 
a function $\phi_n\in \CB'(\mu_1,\mu_2)$ by
$$
\phi_n\biggl(\left(\begin{matrix} a & c\\ 0 & b\end{matrix}\right)k(\theta)\biggr) = 
\mu_1(a)\mu_2(b)(|a/b|)^{1/2}e^{in\theta}
$$
where
$$
k(\theta) = \left(\begin{matrix} \cos\theta & \sin\theta\\ -\sin\theta & \cos\theta\end{matrix}\right)\in K:= O(2)\in G
$$

Let $\CB(\mu_1,\mu_2)\subset \CB'(\mu_1,\mu_2)$ be the (dense) subspace generated by all $\phi_n$. 

Let us describe explicitely the induced action of $\fG = Lie(G)_\BR\otimes\BC$ on $\CB(\mu_1,\mu_2)$. 
Following [L], consider a matrix $A = \left(\begin{matrix} 1 & i\\ 1 & -i\end{matrix}\right)$ so that 
$A^{-1} = \frac{1}{2}\left(\begin{matrix} 1 & 1\\ -i & i\end{matrix}\right)$ (cf. [Ba], (3.5)). 

Let 
$$
X = \left(\begin{matrix} 0 & 1\\ 0 & 0\end{matrix}\right),\ Y = \left(\begin{matrix} 0 & 0\\ 1 & 0\end{matrix}\right),\ 
H = \left(\begin{matrix} 1 & 0\\ 0 & -1\end{matrix}\right) 
$$
Then
$$
A^{-1}XA = \frac{1}{2}\left(\begin{matrix} 1 & -i\\ -i & -1\end{matrix}\right) =: Y',
$$
$$
A^{-1}YA = \frac{1}{2}\left(\begin{matrix} 1 & i\\ i & -1\end{matrix}\right) =: X',
$$
$$
A^{-1}(-iH)A = \frac{1}{2}\left(\begin{matrix} 0 & 1\\ -1 & 0\end{matrix}\right) = k(\pi/2),
$$
or more generally
$$
\left(\begin{matrix} e^{-i\theta} & 0\\ 0 & e^{i\theta}\end{matrix}\right) = Ak(\theta)A^{-1}
$$
Thus if $K' = AKA^{-1}$ then $Lie(K')_\BC = \BC\cdot H$. 

The action of $G$ on $\CB'(\mu_1, \mu_2)$ induces an action of $\fG$ on $\CB(\mu_1, \mu_2)$ which looks as follows:
$$
2X'\phi_n = (s + n + 1)\phi_{n+2},\ 
2Y'\phi_n = (s - n + 1)\phi_{n-2},
\eqno{(1.6.1)} 
$$ 
cf. [JL], Lemma 5.6.

The space $\CB(\mu_1, \mu_2)$ is a {\it $(\fG, K)$-module}, which means that it is a $\fG$-module and 
a $K$-module and the action of $\fk := Lie K$ induced from $\fG$ coincides with the one  
induced from $K$.

\bigskip\bigskip

%\newpage

%\input loke-q-2

\centerline{\bf \S 2. A $q$-deformation.}

\bigskip\bigskip

{\bf 2.1. Category $\CC_q$ and tensor product.} Cf. [Lu]. Let $q$ be a complex number 
different from $0$ and not a root of unity. We fix a value of $\log q$ and for any 
$s\in\BC$ define  $q^s :=  e^{s\log q}$.    
 
Let $U_q = U_q\frak{sl}_2$ denote the 
$\BC$-algebra generated by $E, F, K, K^{-1}$ subject to relations 
$$
K\cdot K^{-1} = 1
$$
$$
KE = q^{2} EK,\ KF = q^{- 2} FK,
$$
$$
[E, F] = \frac{K - K^{-1}}{q - q^{-1}},
$$
cf. [Lu], 3.1.1.

Introduce a comultiplication $\Delta:\ U_q\lra U_q\otimes U_q$ as a unique algebra homomorphism 
such that
$$
\Delta(K) = K\otimes K
$$
$$
\Delta(E) = E\otimes 1 + K\otimes E,\ 
$$
$$
\Delta(F) = F\otimes K^{-1} + 1\otimes F,
$$
cf. [Lu], Lemma 3.1.4.

Let $\CC_q$ denote the category of $\BZ$-graded $U_q$-modules $M = \oplus_{i\in \BZ} M_i$ such that
$$
K x = q^i x,\ x\in M_i.
$$
The comutliplication $\Delta$ above makes $\CC_q$ a tensor category. 

In particular if $M_i$, $i = 1, 2, 3$, are objects of $\CC_q$ then their tensor product 
$M = M_1\otimes M_2\otimes M_3$ is defined; as a vector space it is the tensor product 
of vector spaces underlying $M_i$. The action of $U_q$ is given by
$$
K(x\otimes y\otimes z) = Kx\otimes Ky\otimes Kz 
$$
$$
E(x\otimes y\otimes z) = Ex\otimes y\otimes z + Kx\otimes Ey\otimes z + Kx\otimes Ky\otimes Ez
$$
$$
F(x\otimes y\otimes z) = Fx\otimes K^{-1}y\otimes K^{-1}z + x\otimes Fy\otimes K^{-1}z + x\otimes y\otimes Fz
$$

{\bf 2.2. Infinitesimal principal series.} Set
$$
[n]_q = \frac{q^{n} - q^{-n}}{q - q^{-1}}, 
$$
$q\in\BR_{> 0}, s\in\BC$. 
Thus
$$
\lim_{q\ra 1}\ [n]_q = n.
$$

Let $s\in \BC,\ \epsilon \in\{0, 1\}$. 
Define an object $M_q(s,\epsilon)\in \CC_q$ as follows. As a $\BZ$-graded vector space 
$M_q(s,\epsilon) = \oplus M_i$ where $M_i = \BC\cdot v_i$ if $i\in \epsilon + 2\BZ$ and $0$ otherwise.

An action of the operators $E, F$ are given by 
$$
E v_n = [(s + n + 1)/2]_q v_{n+2},\ F v_n = [(s - n + 1)/2]_q v_{n-2}.
$$
One checks that 
$$
[E_q, F_q] = \frac{q^H - q^{-H}}{q - q^{-1}} = \frac{K - K^{-1}}{q - q^{-1}}
$$
where
$$
Kv_n = q^n v_n,
$$
so $M_q(s,\epsilon)$ is an $U_q$-module. 

%Here is a different way to write down the same representaion. Consider the space of Laurent polynomials 
%$M = \BC[x, x^{-1}]$ and the following two operators acting on it:  
%$$
%Ef(x) = x^2\cdot \frac{q^{(s+1)/2}f(q^{1/2}x) - q^{-(s+1)/2}f(q^{-1/2}x)}{q - q^{-1}}
%$$ 
%and 
%$$
%Ff(x) = x^{-2}\cdot \frac{q^{(s+1)/2}f(q^{-1/2}x) - q^{-(s+1)/2}f(q^{1/2}x)}{q - q^{-1}}
%$$
%They give rise to an $U$-module structure on $M$. These formulas are close to the ones from 
%[PT], [KLS]. 
%After identification $x^n$ with $v_n,\ n\in\BZ$, $M$ gets isomorphic to $M(s,0)\oplus M(s,1)$. 

{\bf 2.2.1. Lemma.} {\it If $s - \epsilon\notin 2\BZ + 1$ then 
$M_q(s,\epsilon)$ is an irreducible $U_q$-module.} 

The proof is the same as in the non-deformed case (see Lemma 1.1.1).

{\bf 2.2.2. The reflection operator.} As in 1.1.2, let us construct an intertwining 
operator
$$
R_q(s):\ M_q(s,\epsilon) \iso M_q(-s,\epsilon).
$$
Suppose that 
$$
R_q(s)v_n = r_n v_n
$$
for some $r_n\in\BC$. As in {\it loc. cit.}, $R_q(s)$ is $U_q$-equivariant iff the numbers 
$r_n$ satisfy the equation 
$$
r_{n+2} = \frac{[(-s+n+1)/2]_q}{[(s+n+1)/2]_q}r_n.
\eqno{(2.2.2.1)}
$$
Suppose we have found a function $\phi(x),\ x\in\BC$, satisfying a functional equation 
$$
\phi(x+1) = [x]_q\phi(x).
\eqno{(2.2.2.2)}
$$
Then
$$
r_n = \frac{\phi((-s+n+1)/2)}{\phi((s+n+1)/2)}
$$
satisfies $(2.2.2.1)$. 

Suppose that $|q| < 1$. In that case consider 
the $q$-Gamma function defined by a convergent infinite product
$$
\Gamma_q(x) = (1 - q)^{1-x}\frac{(q;q)_\infty}{(q^x;q)_\infty}
$$
where 
$$
(a;q)_\infty = \prod_{n=0}^\infty (1 - aq^n),
$$
cf. [GR]. It satisfies the functional equation
$$
\Gamma_q(x+1) = \frac{q^x - 1}{q - 1} \Gamma_q(x).
$$ 
It follows that a function 
$$
\phi(x) = q^{a(x)}\Gamma_{q^2}(x)
$$
satisfies $(2.2.2.2)$ if $a(x)$ satisfies 
$$
a(x+1) - a(x) = 1 - x,
$$
for example 
$$
a(x) = - \frac{x^2}{2} + \frac{3x}{2}
$$
Thus, {\it if we set 
$$
\phi(x) = q^{-(x^2 - 3x)/2}\Gamma_{q^2}(x),
$$
the operator $R_q(s)$ defined by
$$
R_q(s)v_n = \frac{\phi((-s+n+1)/2)}{\phi((s+n+1)/2)}v_n
$$
is an isomorphism $R_q(s):\ M_q(s,\epsilon) \iso M_q(-s,\epsilon)$ 
in $\CC_q$.

It possesses the unitarity property
$$
R_q(-s)R_q(s) = \Id_{M(s)}.
$$}

If $|q| = 1$ then a solution to the functional equation (2.2.2.2) 
may be given in terms of the Shintani-Kurokawa {\it double sine function} (aka Ruijsenaars hyperbolic Gamma function), cf. 
[NU], Prop. 3.3, 
[R], Appendix A. This function is a sort of a "modular double" of $\Gamma_q$.

{\bf 2.3. Theorem} {\it Let $M^i = M_q(s_i,\epsilon_i)\in \CC_q$ be $3$ objects as above such that $s_i - \epsilon_i\notin 
2\BZ + 1$, $i = 1, 2, 3$. 
  
There exists a unique, up to a scalar multiple, function 
$$
f:\ M := M^1\otimes M^2\otimes M^3\lra \BC
\eqno{(2.3.1)}
$$
such that
$$
f(X m) = 0,\ X \in E, F,\ m\in M,
\eqno{(2.3.2)}
$$
$$
f(Km) = f(m),
$$
$$
f(\omega m) = f(m)
\eqno{(2.3.3)}
$$
where $\omega:\ M\iso M$ is an automorphism defined by 
$$
\omega(v_n\otimes v_m\otimes v_k) = v_{-n}\otimes v_{-m}\otimes v_{-k}.
$$}

{\bf 2.4. Proof (beginning).} The argument below is a straightforward generalization of the argument from [L], \S 2.  
The condition $f(Km)=f(m)$ implies that $f(v_n\otimes v_m\otimes v_k)=0$ unless $n+m+k=0$. 
Let us denote 
$$
a_q(n,m) := f(v_n\otimes v_m\otimes v_{-n-m-2}).
$$
The condition $f(\omega m) = f(m)$ gives 
$$
a_q(n,m)=a_q(-n,-m)
\eqno{(2.4.0)}
$$
Since 
$$
f(E(v_n\otimes v_m\otimes v_{-n-m-2})) = 0,
$$
we get
$$
[(s_1+n+1)/2]_q a_q(n+2,m) + 
q^n[(s_2+m+1)/2]_q a_q(n,m+2) +  
$$
$$
+ q^{n+m}[(s_3-n-m-1)/2]_q a_q(n,m) = 0
\eqno{(2.4.1)}
$$
or
$$
[(s_3-n-m-1)/2]_qa_q(n,m) = 
$$
$$
-q^{-n-m}[(s_1+n+1)/2]_q a_q(n+2,m) - q^{-m}[(s_2+m+1)/2]_q a_q(n,m+2)
\eqno{(2.4.1)'}
$$
Similarly, 
$$
f(F(v_n\otimes v_m\otimes v_{-n-m+2}))=0
$$
implies
$$
q^{n-2}[(s_1-n+1)/2]_qa_q(n-2,m)+
$$
$$
q^{n+m-2}[(s_2-m+1)/2]_q a_q(n,m-2)+[(s_3+n+m-1)/2]_qa_q(n,m)=0
\eqno{(2.4.2)}
$$
or
$$
[(s_3+n+m-1)/2]_qa_q(n,m)=
$$
$$
-q^{n-2}[(s_1-n+1)/2]_qa_q(n-2,m)-q^{n+m-2}[(s_2-m+1)/2]_q a_q(n,m-2)
\eqno{(2.4.2)'}
$$
It follows from $(2.4.1)'$:
$$
[(s_3-n-m+1)/2]_qa_q(n-2,m)=
$$
$$-q^{-n-m+2}[(s_1+n-1)/2]_q a_q(n,m)-q^{-m}[(s_2+m+1)/2]_q a_q(n-2,m+2)(5)
\eqno{(2.4.3)}
$$
and
$$
[(s_3-n-m+1)/2]_qa_q(n,m-2)=
$$
$$
-q^{-n-m+2}[(s_1+n+1)/2]_q a_q(n+2,m-2)-q^{-m+2}[(s_2+m-1)/2]_q a_q(n,m)
\eqno{(2.4.4)}
$$
(One could write $(2.4.3) = (2.4.1)'_{n-2,m}$ and $(2.4.4) = (2.4.1)'_{n,m-2}$

Sustitute $(2.4.3)$ and $(2.4.4)$ into $(2.4.2)'$:
$$
\bigg([(s_3-n-m+1)]_q[(s_3+n+m-1)/2]_q-q^{-m}[(s_1+n-1)/2]_q[(s_1-n+1)/2]_q -
$$
$$
-q^n[(s_2-m+1)/2]_q[(s_2+m-1)/2]_q\bigg) a_q(n,m)
$$
$$
=q^{n-m-2}[(s_1-n+1)/2]_q[(s_2+m+1)/2]_q a_q(n-2,m+2)+
$$
$$
[(s_2-m+1)/2]_q[(s_1+n+1)/2]_qa_q(n+2,m-2) 
\eqno{(2.4.5)}
$$
This is a $q$-deformed (1.3.3).

Now comes the main point.

{\bf 2.5. Lemma.} {\it Let $N\in \mathbb Z$ be such that $N\equiv \epsilon_1+\epsilon_2\ (\mod 2)$. 
Suppose we are given $a_q(n,m)$ for $n+m=N$ and they satisfy $(2.4.5)$. Using $(2.4.2)$ let us  
define $a_q(n,m)$ for $n+m=N+2k (k\geq 1)$ inductively. 

Then 
$a_q(n,m)$ satisfies $(2.4.1)$ for $n+m\geq N$.} 

{\bf Proof.} We will prove the lemma by induction on $n+m$. \\
By induction we assume that $(2.4.1)$ is satisfied for all $n+m\leq N-2$. Hence $a_q(n,m)$ also satisfies $(2.4.5)$ for all
$n+m\leq N-2$. 

Let $n+m = N-2$; we want to prove $(2.4.1)$ where $a_q(n+2,m)$ and $a_q(n,m+2)$ are defined from 
$(2.4.2)'$:
$$
ta_q(n+2,m)=-q^n[(s_1-n-1)/2]_qa_q(n,m)-q^{n+m}[(s_2-m+1)/2]_q a_q(n+2,m-2)
\eqno{(2.5.1)}
$$
$$
ta_q(n,m+2)=-q^{n-2}[(s_1-n+1)/2]_qa_q(n-2,m+2)-q^{n+m}[(s_2-m-1)/2]_q a_q(n,m)
\eqno{(2.5.2)}
$$
where $t=[(s_3+n+m+1)/2]_q \neq 0$ by assumption.

We put $(2.5.1)$ and $(2.5.2)$ into the right hand side of $(2.4.1)'$:
$$
-q^{-n-m}[(s_1+n+1)/2]_q a_q(n+2,m)-q^{-m}[(s_2+m+1)/2]_q a_q[n,m+2]=
$$
$$
=-q^{-n-m}[(s_1+n+1)/2]_qt^{-1}
\times
$$
$$
\bigg(-q^n[(s_1-n-1)/2]_qa_q(n,m)-q^{n+m}[(s_2-m+1)/2]_q a_q(n+2,m-2)\bigg)-
$$
$$
-q^{-m}[(s_2+m+1)/2]_qt^{-1}\times
$$
$$
\bigg(-q^{n-2}[(s_1-n+1)/2]_qa_q(n-2,m+2)-q^{n+m}[(s_2-m-1)/2]_q a_q(n,m)\bigg)
$$
$$
=t^{-1}\bigg(q^{-m}[(s_1+n+1)/2]_q[(s_1-n-1)/2]_qa_q(n,m)+
$$
$$
[(s_1+n+1)/2]_q[(s_2-m+1)/2]_q a_q(n+2,m-2)+
$$
$$
q^{n-m-2}[(s_2+m+1)/2]_q[(s_1-n+1)/2]_qa_q(n-2,m+2)+
$$
$$
q^n[(s_2+m+1)/2]_q[(s_2-m-1)/2]_qa_q(n,m)\bigg)
$$
(we substitute (2.4.5) for the second and third terms)
$$
=t^{-1}\bigg(q^{-m}[(s_1+n+1)/2]_q[(s_1-n-1)/2]_q+[(s_3-n-m+1)/2]_q[(s_3+n+m-1)/2]_q-
$$
$$
-q^{-m}[(s_1+n-1)/2]_q[(s_1-n+1)/2]_q-q^n[(s_2-m+1)/2]_q[(s_2+m-1)/2]_q+
$$
$$
q^n[(s_2+m+1)/2]_q[(s_2-m-1)/2]_q\bigg)a_q(n,m)
$$
$$
=t^{-1}\bigg(q^{s_3}+q^{-s_3}-q^{m+n-1}-q^{-m-n+1}-q^{-m+n+1}-q^{-m-n-1}+q^{-m+n-1}+q^{-m-n+1}-
$$
$$
-q^{m+n+1}-q^{-m+n-1}+q^{m+n-1}+q^{-m+n+1}\bigg)a_q(n,m)
$$
$$
=t^{-1}\bigg(q^{s_3}+q^{-s_3}-q^{-m-n-1}-q^{m+n+1}\bigg)a_q(n,m)
$$
$$
=t^{-1}[(s_3+n+m+1)/2]_q[(s_3-n-m-1)/2]_q a_q(n,m)
$$
$$
=[(s_3-n-m-1)/2]_q a_q(n,m)
$$
But this is exactly $(2.4.1)'$! 
This proves  the lemma. $\square$

{\bf 2.6. End of the proof of Thm. 2.3.} By $(2.4.5)$ and equality $a_q(2,-2)= a_q(-2,2)$ we have 
$$
\bigg([(s_3+1)/2]_q[(s_3-1)/2]_q-
$$
$$
[(s_1+1)/2]_q[(s_1-1)/2]_q-[(s_2+1)/2]_q[(s_2-1)/2]_q\bigg)a_q(0,0)=
$$
$$=2[(s_1+1)/2]_q[(s_2+1)/2]_qa_q(2,-2)
\eqno{(2.6.1)}
$$
Let us construct a solution $a_q$ of equations (2.4.0) - (2.4.2) as follows. 

(i) If $\epsilon_1 = 0$, we start from an arbitrary value of $a_q(0,0)$ and define $a_q(2,-2)$ by $(2.6.1)$. 

(ii) If $\epsilon_1 = 1$, we start from an arbitrary value of $a_q(1,-1)$ and set $a_q(-1,1) = a_q(1,-1)$. 

Using $(2.4.5)$, repeatedly, we determine $a_q(n,-n)$ for all positive $n$. Using $(2.4.0)$, we determine $a_q(n,-n)$ for all $n\leq 0$.
Applying $(2.4.2)$ inductively one defines $a_q(n,m)$ for all $n+m>0$. Finally $(2.4.0)$ gives $a_q(n,m)$ for $n+m<0$. 

From the construction, $a_q(n,m)$ satisfies $(2.4.0)$ and $(2.4.2)$ if $n+m>0$ and $(2.4.1)$ if $n+m<0$. Lemma 2.5 shows that $(2.4.1)$ 
is satisfied when $n+m\geq 0$. This proves the existence of $a_q$.  

Since $a_q(n,m)$ is completely determined by its value at $a_q(0,0)$ or $a_q(1,-1)$, the dimension of the space of solutions of the system 
(2.4.0) - (2.4.2) is equal to $1$. This completes the proof of Thm. 2.3. $\square$

\bigskip\bigskip

\centerline{\bf References}

\bigskip\bigskip

[Ba] V.Bargmann, Irreducible unitary representations of the Lorentz 
group, {\it Ann. Math.} {\bf 48} (1947), 568 - 640.

[BR] J.Bernstein, A.Reznikov, Estimates of automorphic functions, 
{\bf 4} (2004), 19 - 37.

[BS] B.V.Binh, V.Schechtman, Remarks on a triple integral, 
arXiv:1204.2117, {\it Moscow Math. J.}, to appear. 

%[vdB] F.J. van de Bult, Ruijsenaars' hypergeometric function and the modular double of $\CU_q(\frak{sl}_2(\BC))$, 
%{\it Adv. Math.} {\bf 204} (2006), 539 - 571. 

[JL] H.Jacquet, R.Langlands, Automorphic forms on $GL(2)$.

[FM] B.Feigin, F.Malikov, Integral intertwining operators and complex powers of differential 
($q$-difference) operators, pp. 15 - 63  in: Unconventional Lie algebras,  
{\it Adv. Soviet Math.} {\bf 17}, Amer. Math. Soc., Providence, RI, 1993.

[GR] G.Gasper, M.Rahman, Basic hypergeometric series.

%[KLS] S.Kharchev, D.Lebedev, M.Semenov-Tian-Shansky, Unitary representations of 
%$\CU_q(\fsl(2,\BR))$, the modular double, and the multiparticle $q$-deformed 
%Toda chains, {\it Comm. Math. Phys.} {\bf 225} (2003), 573 - 609.    

[L] H.Y.Loke, Trilinear forms of $\frak{gl}_2$, {\it Pacific J. Math.}, 
{\bf 197} (2001), 119 - 144. 

[Lu] G.Lusztig, Introduction to quantum groups. 

[NSU] A.Nikiforov, S.Suslov, V.Uvarov, Classical ortothogonal polynomials 
of a discrete variable.

[NU] M.Nishizawa, K.Ueno, Integral solutions of $q$-difference equations 
of hypergeometric type with $|q| = 1$, arXiv:q-alg/9612014. 

[Ok] A.Oksak, Trilinear Lorentz invariant forms, {\it Comm. Math. Phys.} {\bf 29} (1973), 189 - 217. 

%[PT] B.Ponsot, J.Teschner, Clebsch-Gordon and Rakah-Wigner coefficients for a continuos series 
%of representations of $\CU_q(\frak{sl}_2(\BR))$, {\it Comm. Math. Phys.} {\bf 224} (2001), 613 - 655. 

[R] S.Ruijsenaars, A generalized hypergeometric function satisfying four analytic difference equations of Askey-Wilson type, 
{\it Comm. Math. Phys.} {\bf 206} (1999), 639 - 690.  

%[S] J.Stockman, Generalized Cherednik-Macdonald identities, arXiv:0708.0934.

\bigskip\bigskip

Institut de Math\'ematiques de Toulouse, Universit\'e Paul Sabatier, 31062 Toulouse, France

%\newpage

\end{document}